\documentclass{amsart}

\usepackage[english]{babel}
\usepackage{graphicx}
\usepackage[colorlinks]{hyperref}
\usepackage[dvipsnames]{xcolor}

\usepackage[backgroundcolor=white,linecolor=red,bordercolor=red]{todonotes}
\usepackage[shortlabels]{enumitem}
\setenumerate{label=(\roman*)}

\usepackage{hyperref}
\newcommand\myshade{100}
\hypersetup{
  linkcolor  = RoyalBlue!\myshade!black,
  citecolor  = ForestGreen!\myshade!black,
  urlcolor   = RedOrange!\myshade!black,
  colorlinks = true,
}

\usepackage{amsmath}
\usepackage{amssymb}
\usepackage{amsfonts}
\usepackage{amsthm}
\usepackage{mathtools}
\usepackage{bm}
\usepackage{bbold}
\usepackage{tikz-cd}
\usepackage{thmtools}

\usepackage[capitalise]{cleveref}
\crefformat{equation}{(#2#1#3)}
\crefname{subsection}{Section}{Sections}
\crefname{subsubsection}{Section}{Sections}

\theoremstyle{plain}
\newtheorem{thm}{Theorem}[section]
\newtheorem{lem}[thm]{Lemma}
\newtheorem{prop}[thm]{Proposition}
\newtheorem{cor}[thm]{Corollary}

\theoremstyle{definition}
\newtheorem{defn}[thm]{Definition}
\newtheorem{example}[thm]{Example}
\newtheorem{question}[thm]{Question}

\theoremstyle{remark}
\newtheorem{rem}[thm]{Remark}

\makeatletter
\def\@fnsymbol#1{\ensuremath{\ifcase#1\or \dagger\or \ddagger\or
           \dagger\dagger
           \or \ddagger\ddagger \else\@ctrerr\fi}}
\makeatother

\newcommand\N{\ensuremath{\mathbb{N}}}
\newcommand\R{\ensuremath{\mathbb{R}}}

\newcommand{\el}[1]{\ensuremath{\ell_{#1}}}

\DeclareMathOperator{\barespn}{span}

\DeclareMathOperator{\Fin}{Fin}

\usepackage[backend=biber,giveninits=true,style=numeric,sortcites=true,maxnames=100]{biblatex}
\usepackage{csquotes}

\newcommand{\one}{\ensuremath{\mathbb{1}}}
\newcommand{\falg}{$f\!$-algebra}
\newcommand{\lalg}{$\ell$-algebra}
\newcommand{\lsubalg}{$\ell$-subalgebra}

\newcommand{\ev}{\operatorname{ev}}
\newcommand{\cph}[1][K_X]{C_{ph}(#1)}
\newcommand{\czero}[1][K_X^{*}]{C_b^{(0)}(#1)}

\title{$f$-algebra products on AL and AM-spaces}
\author{David Muñoz-Lahoz}
\address{Instituto de Ciencias Matemáticas\\Universidad Autónoma de
Madrid}
\email{david.munnozl (at) uam (dot) es}
\thanks{Research supported by an FPI–UAM 2023 contract (funded by
Universidad Autónoma de Madrid) and by grants PID2020-116398GB-I00 and
CEX2023-001347-S (funded by MCIN/AEI/10.13039/501100011033).}
\date{\today}

\begin{document}

\begin{abstract}

We characterize all \falg\ products on AM-spaces by constructing
a canonical AM-space $W_X$ associated to each AM-space $X$, such that
the \falg\ products on $X$ correspond bijectively to the positive
cone $(W_X)_+$. This generalizes the classical description of
\falg\ products on $C(K)$ spaces. We also identify the unique
product (when it exists) that embeds $X$ as a closed
subalgebra of $C(K)$, and study AM-spaces for which this
product exists---the so-called AM-algebras. Finally, we investigate
AM-spaces that admit only the zero product, providing a
characterization in the AL-space case and examples showing that no
simple characterization exists in general.

\end{abstract}

\maketitle

\section{Introduction}

An \falg\ is a vector lattice $X$ endowed with a real algebra
structure in which $a\wedge b=0$ implies
\[
    (ca)\wedge b=0=(ac)\wedge b\quad\text{for all }c \in X_+.
\]
G.\ Birkhoff and R.\ S.\ Pierce introduced \falg s in
\cite{birkhoff_pierce1956}. This structure, along with the related notion
of $f\!$-rings, has been extensively studied since then (see the surveys and monographs
\cite{johnson1960,huijsmans1991,henriksen1997,boulabiar_buskes_triki2007,depagter1981,keimel1995};
for more recent results on the topic, see
\cite{madden2011,buskes_wickstead2017,buskes_schwanke2018,de_jeu2021}).

A Banach lattice algebra is a Banach lattice together
with a Banach algebra structure in which the product of positive
elements is positive. A.\ W.\ Wickstead asked in \cite[Section
4]{wickstead2017_questions} for a description of the most general
Banach lattice algebra structure on AM-spaces and AL-spaces. This
motivated us to describe the \falg\ products on these spaces.
It will turn out that these products have elegant characterizations; it
is not clear at all whether the general situation will also be so
tractable.

More precisely, it is shown in \cref{thm:fproductsAM} that for every AM-space $X$ there
exists an associated AM-space $W_X$ such that the \falg\ products on
$X$ are in one-to-one correspondence with $(W_X)_+$.
\Cref{sec:fproductsAM} is devoted to the construction of $W_X$ and the
proof of \cref{thm:fproductsAM}. This result generalizes
the well-known fact that the \falg\ products on $C(K)$, where $K$ is a compact Hausdorff space, are in one-to-one correspondence
with $C(K)_+$ (see \cref{prop:falg_prod_C}). There
also exist results in the literature characterizing all the almost \falg\
products \cite{scheffold1981} and all the $d$-algebra
products \cite{boulabiar2004} on $C(K)$.

By a classical theorem of Kakutani \cite[Theorem 1]{kakutani1941}, every AM-space $X$ embeds lattice
isometrically in a space of continuous functions $C(K)$. Since every closed
subalgebra of $C(K)$ is also a vector sublattice, a natural stronger question is
to characterize which products on $X$ can be embedded lattice and algebra
isometrically in a space of continuous functions. By
\cref{thm:charAMalg}, this question turns
out to be equivalent to characterizing when $X$ can be isometrically embedded as a
closed (order and ring) ideal of $C(K)$.

The products on an AM-space that embed in a space of continuous
functions were characterized
intrinsically in \cite{munoz-lahoz_tradacete2025}. In short, these products are the ones for
which $(B_X)_+$ is an approximate algebraic identity, where $B_X$
denotes the unit ball of $X$. It was not
considered in \cite{munoz-lahoz_tradacete2025} how many products
satisfying this property exist. It is shown in
\cref{sec:AMalgebras} that such a product does not always exist but, when it does, it
is unique (\cref{cor:unique}). Following \cite{munoz-lahoz_tradacete2025}, the AM-spaces with this product
are called AM-algebras. In particular, the AM-algebra product on
$C(K)$ is the pointwise product. In \cref{thm:AMalg_prod} it is shown that the AM-algebra product
enjoys nice
characterizations, similar to those of the pointwise product in a
space of continuous functions.

One of the characterizations of an AM-algebra $X$ is that the space
$W_X$ is, in some sense, ``maximal.'' In
\cref{sec:only_zero_prod}, we look at the other end and study those AM-spaces for which $W_X$
is minimal. Note that every vector lattice admits an \falg\ product,
namely, the zero product. Thus we are looking at the vector lattices
that admit a unique \falg\ product. In \cref{sec:only_zero_prod}, several properties of such
spaces are studied. This is motivated by the existence of a simple
characterization of AL-spaces that admit only the zero product: an AL-space admits only the zero product if and only if it
has no atoms (\cref{cor:ALspaces_zero}). This is a simple consequence
of the characterization of the \falg\ products on AL-spaces that is
provided in \cref{thm:ALproducts}. We are able to show that there exist
AM-spaces that only admit the zero product. We are, however, unable to
provide a general and simple characterization of such spaces.

\subsection{Preliminaries}

For the unexplained terminology related to vector and Banach lattices, see
\cite{abramovich_aliprantis2002,aliprantis_burkinshaw2006,meyer-nieberg1991,schaefer1974}.
A \emph{lattice-ordered algebra} (briefly, an \lalg) is an associative
real algebra $X$ that is simultaneously a vector lattice whose
positive cone $X_+$ is closed under multiplication, i.e., $ab \in X_+$
for all $a,b \in X_+$. An \lalg\ $X$ is an \emph{\falg} if
\[
    a\wedge b=0\quad\text{implies}\quad
    (ca)\wedge b=0=(ac)\wedge b\text{ for all }c \in X_+.
\]
We shall also refer to this implication as the \emph{\falg\ property}.
An introduction to the theory of \falg s can be found in \cite[Chapter
20]{zaanen1983}.

A \emph{Banach lattice algebra} is an \lalg\ with a
norm that turns it into a Banach lattice and a Banach algebra. If a
Banach lattice algebra satisfies the \falg\ property, it is
called a \emph{Banach \falg}.

A subspace of an \lalg\ that is at the same time a
vector sublattice and a subalgebra will be referred to as an
\emph{\lsubalg}. Let $X$ and $Y$ be \lalg s.
A map $T\colon X\to Y$ that is at the same time a lattice homomorphism
and an algebra homomorphism will be called an \emph{\lalg\
homomorphism}. Similarly, we define \lalg\ isomorphisms,
isometries and embeddings when $X$ and $Y$ are normed.

Throughout this text, and unless otherwise specified, $K$ will denote
a compact Hausdorff space. Recall that the norm-closed subalgebras of
$C(K)$ are also vector sublattices \cite[Lemma 4.48]{folland1999}, and
that norm-closed order ideals and norm-closed ring ideals coincide:
they are always of the form
\[
    \{\, f \in C(K) : f|_F=0 \, \}
\]
for a certain closed set $F\subseteq K$ (see \cite[Chap.\ III, Sec.\
1, Example 1]{schaefer1974} and \cite[Excercise 70]{folland1999}).

\section{\texorpdfstring{$f$-algebra}{f-algebra} products on AM-spaces}\label{sec:fproductsAM}

The following result characterizing the \falg\ products on AM-spaces
with unit follows from \cite[Theorem 2.2]{conrad1974}.

\begin{prop}\label{prop:falg_prod_C}
    A product $P\colon C(K)\times
    C(K)\to C(K)$ is an \falg\ product if and only if there exists a
    weight $w \in C(K)_+$ such that
    \begin{equation}
    P(f,g)(t)=w(t)f(t)g(t)\quad\text{for all }f,g \in C(K)\text{ and
    }t \in K.
    \end{equation}
    This product is submultiplicative if and only if $\|w\|_\infty \le
    1$.
\end{prop}

In view of the previous proposition, one could say that the \falg\
products on $C(K)$ are parametrized by $C(K)_+$. The goal of
this section is to give a similar characterization for
the \falg\ products on general AM-spaces.

Let $X$ be an AM-space. Denote by $\hom(X,\R)$ the set of functionals
on $X$ that are lattice homomorphisms. The sets
\[
    K_X=\hom(X,\R)\cap B_{X^{*}}, K_X^{*}=K_X\setminus \{0\},\text{ and }K_X^{1}=\hom(X,\R)\cap S_{X^{*}}
\]
are all Hausdorff topological spaces when equipped with the
weak$^*$ topology. Furthermore, $K_X$ is compact. We say that a real
function $f$ defined on $K_X$ is
\emph{positively homogeneous} if $f(\lambda x^{*})=\lambda f(x^{*})$
for every $x^{*} \in K_X$ and $\lambda \in \R_+$ such that $\lambda
x^{*}\in K_X$. For instance, the evaluation $\ev_x$ at some $x \in X$
defined by
\[
\ev_x(x^{*})=x^{*}(x)\text{ for all }x^{*}\in K_X
\]
is a continuous positively homogeneous function on $K_X$.
The set $\cph$ of all continuous positively homogeneous functions on
$K_X$ is a vector lattice with respect to the pointwise addition,
scalar multiplication, and ordering. Furthermore, the norm of uniform
convergence makes $\cph$ an AM-space.

According to a representation result by A.\ Goullet de Rugy
\cite[Corollaire 1.31]{goullet_de_rugy1972}, the map
\[
\begin{array}{cccc}
& X & \longrightarrow & \cph \\
        & x & \longmapsto & \ev_x \\
\end{array}
\]
defines a surjective lattice isometry. From now on, we shall identify
$X$ with $\cph$.

Before characterizing all the \falg\ products on $X$ we need to
compute its centre
\[
\mathcal{Z}(X)=\{\, T \in \mathcal{L}(X) : \pm T\le \lambda I\text{
for some }\lambda >0 \, \}
\]
in terms of $K_X$.

Denote by $\czero$ the space of continuous functions $f\colon
K_X^{*}\to \R$ that are bounded and satisfy $f(\lambda
x^{*})=f(x^{*})$ for all $x^{*}\in K_X^{1}$ and $0<\lambda \le 1$.
Note that this is a closed \lsubalg\  of $C_b(K_X^{*})$. In
particular, it is a unital AM-space with the constant function one as
unit.

For every $h \in \czero$, define the operator $M_h\colon \cph \to
\cph$ by
\[
    M_h(f)(x^{*})=
    \begin{cases}
        h(x^{*})f(x^{*})  &\text{if }x^{*}\in K_X^{*}\\
        0 &\text{if }x^{*}=0
    \end{cases}.
\]
Note that $M_h(f)$ is certainly continuous outside $0$ and that, if
$x_\alpha ^{*}\to 0$ in $K_X$ with $x_\alpha ^{*}\neq 0$ for all
$\alpha $ big enough, then
\[
    M_h(f)(x_\alpha ^{*})=h(x_\alpha ^{*})f(x_\alpha ^{*})\to 0
\]
since $h(x_\alpha ^{*})$ is bounded and $f(x_\alpha ^{*})\to 0$.
This shows that $M_h(f)$ is continuous at $0$.
Also, $M_h(f)$ is a positively
homogeneous function. Hence, the operator $M_h$ is well-defined. Note
that it is also central, since it is linear and satisfies $\pm M_hf\le
\|h\|_\infty f$ for all $f\ge 0$. (Side note: every positively
homogeneous function must vanish at $0$; from now on, we will
not write the value at $0$ explicitly.)

Next we show that every operator in the centre of $\cph$ is of the
form described above. We don't claim any originality in this result:
it is well-known that central operators in $\cph$ are multiplication
operators (see \cite{wils1971,buck1961}). Here, we provide a
characterization suited to our purposes using standard techniques.

\begin{prop}\label{prop:center}
    Let $X$ be an AM-space. The map
    \[
    \begin{array}{cccc}
    M\colon& \czero & \longrightarrow & \mathcal{Z}(\cph) \\
            & h & \longmapsto & M_h \\
    \end{array}
    \]
    is an \lalg\ isometry.
\end{prop}
\begin{proof}
    It is direct to check that $M$ is linear and that it is an algebra
    homomorphism with $M\one_{K_X^{*}}=I$. That it is a lattice
    homomorphism follows from the fact that $|T|f=|Tf|$ for every $T
    \in \mathcal{Z}(\cph)$ and $f \in \cph_+$ (see \cite[Theorem
    3.30]{abramovich_aliprantis2002}). The norm in
    $\mathcal{Z}(\cph)$ is given by $\|T\|=\inf \{\, \lambda >0 : \pm T\le \lambda
    I\, \} $ for all $T \in \mathcal{Z}(\cph)$ (see \cite[Theorem
    3.31]{abramovich_aliprantis2002}). Therefore,
    once we establish that $M$ is bijective, it will follow that it is
    an isometry.

    That $M$ is injective is straightforward.
    To show that it is surjective, we will use an adaptation of the
    arguments in \cite[Section 5]{buck1961} (see also
    \cite[Lemma 2.1]{wickstead1988}). Let $T \in
    \mathcal{Z}(\cph)$ be a central operator. Let $x^{*}\in K_X^{*}$
    and choose $f,g \in \cph$ with $f(x^{*})g(x^{*})\neq 0$. Setting
    \[
    k=f-\frac{f(x^{*})}{g(x^{*})}g\text{ in }\cph
    \]
    we get, as $T$ is central,
    \[
    |(Tk)(x^{*})|\le \|T\| |k(x^{*})|=0.
    \]
    Whence,
    \[
    \frac{(Tf)(x^{*})}{f(x^{*})}=\frac{(Tg)(x^{*})}{g(x^{*})}.
    \]
    Accordingly, the function $h\colon K_X^{*}\to \R$ given by
    \[
    h(x^{*})=\frac{(Tf)(x^{*})}{f(x^{*})}\text{ for all }x^{*} \in
    K_X^{*},
    \]
    where $f$ is chosen arbitrarily in $\cph$ with $f(x^{*})\neq 0$, is
    well-defined. Since $f \in \cph$ is continuous, if $f(x^{*})\neq 0$
    there exists a weak$^*$ open neighbourhood $U$ of $x^{*}$ in which
    $f|_U$ does not vanish. Then $h|_U=(Tf)|_U/f|_U$ and, since both
    $f$ and $Tf$ are continuous, $h|_U$ is continuous. Hence $h$ is
    locally continuous, and therefore continuous. Also note that
    $h$ is bounded, since
    \[
    |h(x^{*})|=\frac{|Tf(x^{*})|}{|f(x^{*})|}\le
    \|T\|\frac{|f(x^{*})|}{|f(x^{*})|}=\|T\|.
    \]

    Thus far we have constructed $h \in C_b(K_X^{*})$ such that
    $Tf(x^{*})=h(x^{*})f(x^{*})$ for all $x^{*}\in K_X^{*}$ and $f \in
    \cph$ with $f(x^{*})\neq 0$. Note that, if $f(x^{*})=0$, then also
    $Tf(x^{*})=0$ (because $T$ is central), so the identity
    $Tf(x^{*})=h(x^{*})f(x^{*})$ still holds in this case. For an
    element $x^{*}\in K_X^{1}$ and $0<\lambda \le 1$, using the
    homogeneity of $f$ and $Tf$:
    \[
    \lambda h(x^{*})f(x^{*})=\lambda Tf(x^{*})=Tf(\lambda x^{*})=h(\lambda x^{*})f(\lambda
    x^{*})=\lambda h(\lambda x^{*})f(x^{*}).
    \]
    Choosing an $f$ for which $f(x^{*})\neq 0$ shows that
    $h(x^{*})=h(\lambda x^{*})$. Hence $h \in \czero$ and $T=M_h$, as
    wanted.
\end{proof}

To each function $w\colon K_X^{1}\to \R$ associate a map $w_{-1}\colon
K_X^{*}\to \R$ defined by
\[
w_{-1}(x^{*})=\frac{w(x^{*}/\|x^{*}\|)}{\|x^{*}\|},\quad x^{*}\in
K_X^{*}.
\]
Consider
\[
W_X=\{\, w \in C_b(K_X^{1}) : w_{-1}\text{ is continuous on }K_X^{*} \, \}.
\]
We are going to show that $W_X$ is a closed vector sublattice of
$C_b(K_X^{1})$, and therefore an AM-space. The positive cone of $W_X$
will end up being in one-to-one correspondence with the \falg\
products on $X$, i.e., $(W_X)_+$ will be our space of weights. Note
that we are not asking for $w_{-1}$ to be bounded on $K_X^{*}$ (in
many cases it is not; see \cref{sec:AMalgebras}), we are only asking
for it to be continuous.

From $(\lambda w+w')_{-1}=\lambda w_{-1}+(w')_{-1}$
and $|w|_{-1}=|w_{-1}|$ it follows that $W_X$ is a vector sublattice of
$C_b(K_X^{1})$. To show that it is closed, let $(w_n)\subseteq W_X$ be
a sequence converging to $w$ in $C_b(K_X^{1})$. For $\delta >0$, denote
\[
K_X^{\delta }=\{\, x^{*}\in K_X : \|x^{*}\|\ge \delta  \, \}.
\]
Estimate
\begin{align*}
    \sup_{x^{*}\in K_X^{\delta }}
    |(w_n)_{-1}(x^{*})-w_{-1}(x^{*})|&=\sup_{x^{*}\in K_X^{\delta
    }}\frac{|w_n(x^{*}/\|x^{*}\|)-w(x^{*}/\|x^{*}\|)|}{\|x^{*}\|}\\
    &\le \frac{\|w_n-w\|_{C_b(K_X^{1})}}{\delta }.
\end{align*}
This shows that $(w_n)_{-1}$ converges uniformly to $w_{-1}$ on
$K_X^{\delta }$. In particular, $w_{-1}$ is continuous on $K_X^{\delta
}$. Since this is true for every $\delta >0$, it follows that
$w_{-1}$ is continuous on $K_X^{*}$. Hence $w \in W_X$ and $W_X$ is a
closed vector sublattice of $C_b(K_X^{1})$.

Every element of $(W_X)_+$ defines an \falg\ product on
$X$ in the following way:

\begin{lem}
    Let $X$ be an AM-space and let $w \in (W_X)_+$. Put
    \[
    \begin{array}{cccc}
    P\colon& \cph \times \cph & \longrightarrow & \cph \\
           & (f,g) & \longmapsto & \chi _{\{0\}}w_{-1}fg \\
    \end{array},
    \]
    where $\chi _{\{0\}}$ is the characteristic function of $\{0\}$.
    Then $P$ is an
    \falg\ product on $X$. If $\|w\|_\infty \le 1$, then $P$ is
    submultiplicative.
\end{lem}
\begin{proof}
    The only non-trivial point is to check that, given $f,g \in \cph$,
    the function $P(f,g)$ is continuous at $0$. But this follows
    readily from the inequalities
    \[
    |P(f,g)(x^{*})|=\left| w\left( \frac{x^{*}}{\|x^{*}\|} \right)
    f\left( \frac{x^{*}}{\|x^{*}\|} \right) g(x^{*}) \right| \le
    \|w\|_\infty \|f\| |g(x^{*})|
    \]
    that hold for every $x^{*}\in K_X^{*}$. Notice here that the same
    inequalities yield straightforwardly the submultiplicativity
    result.
\end{proof}

The main result of this paper is that every \falg\ product on $X$ is
of this form.

\begin{thm}\label{thm:fproductsAM}
    Let $X$ be an AM-space identified with $\cph$. A product
    \[
        P\colon \cph \times \cph \to \cph
    \]
    is an \falg\ product if and only if there
    exists a weight $w\in  (W_X)_+$ such that
    \[
        P(f,g)=\chi _{\{0\}}w_{-1}fg
    \]
    for all $f,g \in \cph$. This product is submultiplicative if and only if $\|w\|_\infty \le 1$.
\end{thm}
\begin{proof}
    After previous lemma, it only remains to show that
    every \falg\ product is of the desired form and that, if the
    product is submultiplicative, then $\|w\|_\infty \le 1$.
    Let $P\colon \cph \times \cph \to \cph$ be an \falg\ product. For
    every $g \in \cph$, the \falg\ property implies that the map
    $P(.,g)$ is an orthomorphism. By
    \cite[Theorem 3.29]{abramovich_aliprantis2002},
    $P(.,g)\in \mathcal{Z}(\cph)$, and by
    \cref{prop:center} there exists a unique function $h_g \in
    \czero$ such that $P(.,g)=M(h_g)$. Consider the map
    \[
    \begin{array}{cccc}
    H\colon& \cph & \longrightarrow & \czero \\
            & g & \longmapsto & h_g \\
    \end{array}.
    \]
    That $H$ is a lattice homomorphism follows from a routine
    calculation involving the injectivity of $M$ and the fact that
    positive central operators are lattice homomorphisms.

    Let $x ^{*} \in K_X^{*}$ and let $\delta _{x^{*}} \in
    \czero^{*}$ denote the usual evaluation at $x^{*}$
    functional. The map $H^{*}\delta _{x^{*}}=\delta
    _{x^{*}}\circ H\colon
    \cph\to \R$ is a lattice homomorphism. Identifying $X$ with
    $\cph$ (through $x\mapsto \ev_x$), the definition of $K_X^{*}$ implies that there exist
    $\phi (x^{*}) \in K_X^{*}$ and $w_1(x^{*})\in \R_+$ such that $H^{*}\delta
    _{x^{*}}=w_1(x^{*})\phi
    (x^{*})$. This means that, for every $x \in X$, $(H^{*}\delta
    _{x^{*}})(\ev_x)=w_1(x^{*})\phi (x^{*})(x)$. A mere rewrite
    yields
    \[
        (H \ev_x)(x^{*}) = w_1(x^{*}) \ev_x( \phi (x^{*})).
    \]
    This way we can define maps $\phi \colon K_X^{*}\to K_X$
    and $w_1\colon K_X^{*}\to \R_+$ in the obvious way
    satisfying $Hg=w_1 g\circ \phi $ on $K_X^{*}$ for every $g \in
    \cph$. The
    positive homogeneity of the functions implies that $P$ has the form:
    \begin{equation}\label{eq:falg_prod}
    P(f,g)(x^{*})=w_1(x^{*})f(x^{*})g(\phi (x^{*}))\quad\text{for all
    }x^{*}\in K_X^{*}.
    \end{equation}

    Next we want to show that, for every $x^{*} \in K_X^{*}$ for which
    $w_1(x^{*})\neq 0$, $\phi
    (x^{*})$ is a scalar multiple of $x^{*}$. Otherwise, $x^{*}$ and
    $\phi (x^{*})$ must be disjoint in $X^{*}$, because they are
    non-zero
    lattice homomorphisms, and therefore atoms in $X^{*}$ (see
    \cite[Problem 2.3.6]{abramovich_aliprantis2002problems}). Accordingly, there exist $x,y \in X_+$ such that
    \[
    \phi (x^{*})(x) < x^{*}(x)\quad\text{and}\quad x^{*}(y)<\phi
    (x^{*})(y).
    \]
    Then, using that $w_1(x^{*})>0$,
    \[
    P(\ev_x,\ev_y)(x^{*})=w_1(x^{*})x^{*}(x) \phi (x^{*})(y)>w_1(x^{*})\phi
    (x^{*})(x) x^{*}(y) =P(\ev_y,\ev_x)(x^{*}).
    \]
    This contradicts the fact that every \falg\ product on an
    Archimedean vector lattice is commutative \cite[Theorem
    140.10]{zaanen1983}. Hence $\phi (x^{*})$
    must be a positive multiple of $x^{*}$ whenever $w_1(x^{*})\neq
    0$.

    Define
    the map $w_2\colon K_X^{*}\to [0,1]$ by $\phi
    (x^{*})=w_2(x^{*})x^{*}/\|x^{*}\|$,
    whenever $w_1(x^{*})\neq 0$, and by $w_2(x^{*})=0$ (or any other value)
    whenever $w_1(x^{*})=0$. Set $w=w_1w_2$. Then \eqref{eq:falg_prod}
    becomes, using the homogeneity of $g$,
    \[
    P(f,g)(x^{*})=\frac{w(x^{*})}{\|x^{*}\|}f(x^{*})g(x^{*}).
    \]
    By the positive homogeneity of $f$ and $g$,
    \[
        P(f,g)(x^{*}/\|x^{*}\|)=\frac{w(x^{*}/\|x^{*}\|)}{\|x^{*}\|^2}f(x^{*})g(x^{*}).
    \]
    Also
    \[
    \frac{1}{\|x^{*}\|}P(f,g)(x^{*})=\frac{w(x^{*})}{\|x^{*}\|^2}f(x^{*})g(x^{*})
    \]
    and since $P(f,g)$ is positively homogeneous, it follows by
    choosing appropriate $f$ and $g$ that
    $w(x^{*})=w(x^{*}/\|x^{*}\|)$. Hence we can restrict the
    definition of $w$ to $K_X^{1}$ and write
    \[
    P(f,g)(x^{*})=\frac{w(x^{*}/\|x^{*}\|)}{\|x^{*}\|}f(x^{*})g(x^{*})\quad\text{for
    all }x^{*}\in K_X^{*}.
    \]

    It remains to check that $w \in (W_X)_+$. By construction, $w$ is
    positive. Since $P$ is
    positive on both coordinates, it must be bounded by some constant
    $\|P\|$. Consequently, for every $x^{*}\in K_X^{1}$
    and $f \in \cph$ with $\|f\|_\infty \le 1$, we have
    \[
    |w(x^{*})| f(x^{*})^2=|P(f,f)(x^{*})|\le \|P(f,f)\|\le \|P\|.
    \]
    Fixed $x^{*}$, $f(x^{*})$ gets arbitrarily close to $1$ as $f$
    ranges through $\cph_+$, so $|w(x^{*})|\le
    \|P\|$. In particular, if $P$ is
    submultiplicative, then $\|w\|_\infty \le \|P\|\le 1$.

    To show that $w_{-1}$ is continuous, let $x_\alpha ^{*}\to x^{*}$
    be a convergent net in $K_X^{*}$. Let $f \in \cph$ be such that
    $f(x^{*})\neq 0$. Since both $f$ and $P(f,f)$ are continuous on
    $K_X$, we have
    \[
    f(x_\alpha ^{*})^2\to f(x^{*})^2\quad\text{and}\quad w_{-1}(x_\alpha
    ^{*})f(x_{\alpha }^{*})^2\to w_{-1}(x^{*})f(x^{*})^2.
    \]
    Choose $\alpha _0$ such that $|f(x_\alpha ^{*})^2-f(x^{*})^2|\le
    f(x^{*})^2/2$ for all $\alpha \ge \alpha _0$. Then
    \begin{align*}
        |w_{-1}(x_\alpha ^{*})-w_{-1}(x^{*})|&\le \frac{2}{f(x^{*})^2}
        |f(x_\alpha ^{*})|^2|w_{-1}(x_\alpha
        ^{*})-w_{-1}(x^{*})|\\
       &=\frac{2}{f(x^{*})^2}|f(x_\alpha ^{*})^2w_{-1}(x_\alpha
       ^{*})-f(x_\alpha ^{*})^2w_{-1}(x^{*})|\\
       &\le \frac{2}{f(x^{*})^2}|f(x_\alpha
       ^{*})^2w_{-1}(x_\alpha ^{*})-f(x^{*})^2w_{-1}(x^{*})|\\
       &\phantom{\le}+\frac{2}{f(x^{*})^2}|f(x^{*})^2 -f(x_\alpha ^{*})^2||w_{-1}(x^{*})|
    \end{align*}
    for all $\alpha \ge \alpha _0$. Taking limits in $\alpha $ shows
    that $w_{-1}(x_\alpha ^{*})\to w_{-1}(x^{*})$.
\end{proof}

\begin{rem}
    Let $X$ be an AM-space with unit, say $X=C(K)$ for a certain
    compact Hausdorff space $K$.
    Then $K_X^{1}$ is homeomorphic to $K$ (in particular, it is
    compact) and therefore $W_X\subseteq C_b(K_X^{1})=C(K)$. Since the
    norm is continuous on $K_X$ (because every positive functional
    attains its norm at $\one_K$) it follows that $W_X=C(K)$. Hence,
    for AM-spaces with unit, previous theorem reduces to
    \cref{prop:falg_prod_C}.
\end{rem}

So, for AM-spaces with unit, $W_X=C_b(K_X^{1})$, that is, $W_X$ is
maximal. The AM-spaces with the property that $W_X$ is maximal are
surprisingly well-behaved. They are the object of study of the
following section.

\section{AM-algebras}\label{sec:AMalgebras}

According to Kakutani's theorem \cite[Theorem 1]{kakutani1941} every AM-space embeds lattice isometrically in
a space of continuous functions. It is natural to ask: which of the
products in \cref{thm:fproductsAM} embed \lalg\
isometrically in a space of continuous functions? (Note that every
\lsubalg\  of a space of continuous functions is an \falg, so it
suffices to consider \falg\ products to characterize \emph{all}
products on an AM-space that embed in a space of continuous
functions.) The products with this property were characterized
intrinsically in \cite{munoz-lahoz_tradacete2025} using
the following notion.

\begin{defn}
    Let $X$ be a Banach lattice and let $(e_\alpha )\subseteq X_+$ be
    a net. We say that $X$ is an \emph{AM-space with approximate unit}
    $(e_\alpha )$ if $x^{*}(e_\alpha )\to \|x^{*}\|$ for all $x^{*}\in
    X_+^{*}$. We also say that $(e_\alpha )$ is an \emph{approximate order
    unit} of $X$.
\end{defn}

It turns out that having an approximate order unit is equivalent to
being an AM-space \cite[Lemma 2.5]{munoz-lahoz_tradacete2025}. Providing an approximate order unit
is, morally, a quantitative way of saying that a Banach lattice is an
AM-space. This approach to AM-spaces enables us to go one step
further and establish a relation between the AM-space and algebraic
structures. Recall
that a net $(e_\alpha )$ in a Banach algebra $A$ is an \emph{approximate
(algebraic) identity} if both $e_\alpha a\to a$ and $a e_\alpha \to a$ hold for
every $a \in A$.

\begin{defn}
    Let $A$ be a Banach lattice algebra and let $(e_\alpha )\subseteq
    A_+$ be a net. We say that $A$ is an \emph{AM-algebra with
    approximate unit} $(e_\alpha )$ if $(e_\alpha )$ is both an
    approximate order unit and an approximate algebraic identity.

    If $A$ has an order unit $e$ that is at the same time an algebraic
    identity, then $A$ is said to be an \emph{AM-algebra with unit}
    $e$.
\end{defn}

This property characterizes the products on AM-spaces that can be
algebra and lattice isometrically embedded in spaces of continuous functions.

\begin{thm}[{\cite[Theorem 2.10]{munoz-lahoz_tradacete2025}}]
    A Banach lattice algebra $A$ is an AM-algebra with approximate
    unit if and only if it is \lalg\ isometric to a closed
    \lsubalg\  of $C(K)$, for some compact
    Hausdorff space $K$.
\end{thm}

Moreover, we have the following analogue of Kakutani's theorem.

\begin{thm}[{\cite[Theorem 1.2]{munoz-lahoz_tradacete2025}}]
    \label{thm:generalAMalg}
    \begin{enumerate}
         \item[]
        \item Every AM-algebra with unit is \lalg\
            isometric to $C(K)$ for some compact Hausdorff space $K$,
            with the unit corresponding to the constant one function $\mathbb{1}_K$.
        \item Every AM-algebra with approximate unit is
            \lalg\ isometric to a closed \lsubalg\
            of $C(K)$ for some compact Hausdorff space $K$. Moreover,
            choosing $K$ appropriately, there exists a closed set
            $F\subseteq K$ such that the AM-algebra is lattice and
            algebra isometric to the closed \lsubalg\  of
            $C(K)$:
            \[ \{\, f \in C(K) : f(t)=0\text{ for all }t \in
            F\, \}. \]
            In particular, it embeds as an order and
            algebraic ideal in $C(K)$.
    \end{enumerate}
\end{thm}

Even though it was not clear from the definition, a corollary of this
theorem is that, if $A$ is an AM-algebra with an approximate unit,
then $(B_A)_+$ is an approximate identity. Indeed, if we identify $A$
with a closed ideal of $C(K)$, it can be shown that $(B_A)_+$ is an approximate
identity using the standard argument from C$^{*}$-algebras (see, for
instance, \cite[Section II.4]{blackadar2006}).

With this language, the question of which products on an AM-space embed in $C(K)$
becomes: which of the submultiplicative products in \cref{thm:fproductsAM} give the
AM-space $X$ a structure of AM-algebra with approximate unit? Since
$X$ is already an AM-space, these products are precisely the ones
making $(B_X)_+$ into an approximate algebraic identity. The following
general result shows that there can be at most one such product.

\begin{prop}\label{prop:unique_general}
    Let $X$ be a Banach lattice, and let $(e_\alpha )$ be a net in
    $X$. Then there exists at most one Banach \falg\ product on $X$ for
    which $(e_\alpha )$ is an approximate algebraic identity.
\end{prop}
\begin{proof}
    Let $\ast$ and $\star$ be two Banach \falg\ products on $X$.
    Using that central operators commute, it is not difficult to check
    that
    \[
        (f\ast g) \star h = (f\ast h) \star g
    \]
    holds for every $f,g,h \in X$. Suppose $(e_\alpha )$ is an
    approximate algebraic identity for both $\ast$ and $\star$. Then
    \[
        (f\ast g)\star e_\alpha =(f\ast e_\alpha )\star g.
    \]
    Passing to the limit, this shows that $\ast$ and $\star$ are
    equal.
\end{proof}

\begin{cor}\label{cor:unique}
    Let $X$ be an AM-space. If there exists a product making $X$ into an
    AM-algebra with approximate unit, then it is unique.
\end{cor}

When such a product exists, we will
simply say that $X$ is an \emph{AM-algebra} and call the associated
product the \emph{AM-algebra product}. Of course, every AM-space
with unit is an AM-algebra, with the AM-algebra product being the
pointwise product. More generally, the following theorem characterizes the AM-spaces
that are AM-algebras in several ways.

\begin{thm}\label{thm:charAMalg}
    Let $X$ be an AM-space. The following are equivalent:
    \begin{enumerate}
        \item $X$ is an AM-algebra.
        \item $\one_{K_X^{1}} \in W_X$.
        \item The norm is weak$^{*}$ continuous on $K_X^{*}$.
        \item $W_X$ is maximal (i.e., $W_X=C_b(K_X^{1})$).
        \item $X$ is lattice isometric to a closed subalgebra of
            $C(K)$, for a certain compact Hausdorff space $K$.
        \item $X$ is lattice isometric to a closed ideal of $C(K)$,
            for a certain compact Hausdorff space $K$.
        \item $X$ is lattice isometric to $C_0(\Omega )$, for a
            certain locally compact Hausdorff space $\Omega $.
    \end{enumerate}
\end{thm}
\begin{proof}
    Suppose $X$ is an AM-algebra, and let $P$ denote its product. Then
    $P$ is given, in the terminology of \cref{thm:fproductsAM}, by a
    weight $w \in (W_X)_+$ with norm $\|w\|_\infty \le 1$. The unit
    ball $(f_\alpha )=(B_X)_+$ is an approximage algebraic identity
    for $P$,
    that is, $P(f_\alpha,g)\to g$ for all $g \in \cph$. This means in particular that
    \[
        0=\lim_\alpha \sup_{y^{*}\in K_X^{1}} |g(y^{*})| | w(y^{*})
        f_\alpha (y^{*})-1|.
    \]
    Hence $w(y^{*})f_\alpha (y^{*})\to 1$ for all $y^{*} \in K_X^{1}$.
    Since $0\le f_\alpha (y^{*})\le 1$ for all $\alpha $, $w(y^{*})\ge 1$. On
    the other hand, since the product is submultiplicative, $w(y^{*})\le
    1$. Thus $w=\one_{K_X^{1}} \in W_X$, and this shows that \textit{(i)}
    implies \textit{(ii)}.

    If $w=\one_{K_X^{1}} \in W_X$,
    then $w_{-1}(x^{*})=1/\|x^{*}\|$ is continuous on $K_X^{*}$, so
    \textit{(ii)} implies \textit{(iii)}. If the norm is continuous on
    $K_X^{*}$, then for every $w \in C_b(K_X^{1})$ the function
    \[
    w_{-1}(x^{*})=\frac{w(x^{*}/\|x^{*}\|)}{\|x^{*}\|}
    \]
    is continuous on $K_X^{*}$. This shows that \textit{(iii)} implies
    \textit{(iv)}.

    To prove that \textit{(iv)} implies \textit{(v)}, suppose
    $W_X=C_b(K_X^{1})$. In particular, $\one_{K_X^{1}} \in W_X$, and so
    \[
    P(f,g)(x^{*})=\frac{1}{\|x^{*}\|}f(x^{*})g(x^{*})
    \]
    defines an \falg\ product on $\cph$. Consider the restriction map
    \[
    \begin{array}{cccc}
    R\colon& \cph & \longrightarrow & C_b(K_X^{1}) \\
            & f & \longmapsto & f|_{K_X^{1}} \\
    \end{array}.
    \]
    Note that this is a lattice isometric embedding. Moreover, it is
    an algebra homomorphism for the product $P$ in $\cph$ and the
    pointwise product in $C_b(K_X^{1})$. Since $C_b(K_X^{1})$ is
    \lalg\ isometric to $C(\beta K_X^{1})$, where $\beta
    K_X^{1}$ denotes the Stone--Čech compactification of $K_X^{1}$, the desired
    implication follows.

    If $X$ is isometric to a closed
    subalgebra of $C(K)$, then with the product it inherits from
    $C(K)$ it is an AM-algebra with approximate unit, and therefore it is
    isometric to a closed ideal of a space of continuous functions
    (this is part of \cref{thm:generalAMalg}). Hence \textit{(v)}
    implies \textit{(vi)}. Suppose $X$ is lattice isometric to $J_F=\{\, f \in C(K) : f|_F=0 \, \} $ for some closed
    $F\subseteq K$. It is then standard to check that the restriction
    map
    \[
    \begin{array}{cccc}
            & J_F& \longrightarrow & C_0(K\setminus F) \\
            & f & \longmapsto & f|_{K\setminus F} \\
    \end{array}
    \]
    is a (surjective) lattice isometry. This shows that \textit{(vi)}
    implies \textit{(vii)}. That \textit{(vii)} implies \textit{(i)}
    is immediate.
\end{proof}

\begin{rem}\label{rem:pointwise}
    In items \textit{(v)}, \textit{(vi)}, and \textit{(vii)}, the
    AM-algebra product on $X$ is precisely the pointwise product.
    Indeed, in every case the pointwise product makes $X$ into an
    \lalg\ that is lattice isometric to a closed \lsubalg\ of a certain
    $C(K)$. By uniqueness, this product is the
    AM-algebra product on $X$.
\end{rem}

\begin{rem}\label{rem:AMalg_formula}
    If $X$ is an AM-algebra, then the AM-algebra product is given by
    \[
        P(f,g)(x^{*})=\frac{1}{\|x^{*}\|}f(x^{*})g(x^{*})\quad\text{for
        all }f,g \in \cph.
    \]
    In particular, whenever this product is well-defined, the unit
    ball is an approximate algebraic identity for it.
\end{rem}

The following is an example of an AM-space that does not satisfy the
conditions of previous theorem.

\begin{example}
    Consider $X=c$ with the usual linear and lattice structures, but
    with norm
    \[
    \|x\|=(2\lim_n |x_n|)\vee \|x\|_\infty,
    \]
    \todo{it is immediate to check that this is an AM-space}
    where $\|x\|_\infty =\sup_n |x_n|$. Consider the point-evaluation
    functionals $\delta _j$ at coordinates $j \in \N$ and the limit
    functional $\delta _\infty (x)=\lim_n x_n$. These are all
    lattice homomorphisms. It is easy to check that $\|\delta
    _j\|=1$ for every $j \in \N$, whereas $\|\delta _\infty \|=1/2$.
    We claim that $w=\one_{K_X^{1}}\not\in W_X$. Indeed,
    $\delta _j \to \delta _\infty $ in $K_X^{*}$, yet $w_{-1}(\delta
    _j)=1$ and
    \[
    w_{-1}(\delta _\infty )=\frac{w(2 \delta _\infty )}{\|\delta
    _\infty \|}=2.
    \]
    Thus $w_{-1}$ is not continuous on $K_X^{*}$.
\end{example}

The following result distinguishes the AM-algebra product among all
products in several ways.

\begin{thm}\label{thm:AMalg_prod}
    Let $X$ be an AM-algebra and let $P$ be a product on $X$. The following are equivalent:
    \begin{enumerate}
        \item $P$ is the AM-algebra product on $X$.
        \item $(B_X)_+$ is an approximate algebraic identity
            for $P$.
        \item When $X$ is identified with $\cph$:
            \[
            P(f,g)(x^{*})=\frac{1}{\|x^{*}\|}f(x^{*})g(x^{*})
            \]
            for all $f,g \in \cph$ and $x^{*}\in K_X^{*}$.
        \item For every $z^{*}\in K_X^{1}$ and $x,y \in X$:
            \[
                z^{*}(P(x,y))=z^{*}(x)z^{*}(y).
            \]
    \end{enumerate}
    In the case that $X$ is an AM-space with unit $e$, all of the
    above are also equivalent to each of the following additional
    conditions.
    \begin{enumerate}[resume]
        \item $e$ is the algebraic identity of $P$.
        \item When $X$ is identified with a certain $C(K)$, $P$ is the pointwise product.
    \end{enumerate}
\end{thm}
\begin{proof}
    That \textit{(i)} and \textit{(ii)} are equivalent is the main
    content of \cref{thm:generalAMalg}. That \textit{(i)} implies
    \textit{(iii)} is \cref{rem:AMalg_formula}. If \textit{(iii)} holds,
    then
    \[
    P(f,g)(z^{*})=f(z^{*})g(z^{*})
    \]
    for every $z^{*}\in K_X^{1}$; this is precisely \textit{(iv)}.

    Suppose \textit{(iv)} holds. Identify $X$ with $\cph$, so that the
    condition on the product becomes, as before, $P(f,g)(z^{*})=f(z^{*})g(z^{*})$
    for every $z^{*}\in K_X^{1}$ and $f,g \in \cph$. The restriction map
    \[
    \begin{array}{cccc}
    R\colon& \cph & \longrightarrow & C_b(K_X^{1}) \\
            & f & \longmapsto & f|_{K_X^{1}} \\
    \end{array}
    \]
    is a lattice isometric embedding. The condition on the product
    implies that it is also an algebra homomorphism. Since
    $C_b(K_X^{1})$ is \lalg\ isometric to $C(\beta K_X^{1})$,
    \textit{(i)} follows.

    Suppose that $X$ is an AM-space with unit $e$. Suppose $P$ is
    the AM-algebra product on $X$. Since $(B_X)_+$ converges in norm
    to $e$, it follows that $e$ is the algebraic identity of $P$. So
    \textit{(i)} implies \textit{(v)}. If $e$ is the algebraic
    identity of $P$, then $X$ is an AM-algebra with unit $e$, and
    the first item in \cref{thm:generalAMalg} implies \textit{(vi)}.
    That \textit{(vi)} implies
    \textit{(i)} is direct (see \cref{rem:pointwise}).
\end{proof}

The previous two theorems complement very nicely the results of
\cite{munoz-lahoz_tradacete2025}, because they show precisely when an
AM-space can be made into an AM-algebra and, in that case, what the
AM-algebra product looks like. They also provide
practical ways to compute the AM-algebra product in concrete
AM-spaces. The following example illustrates how.

\begin{example}
    Let $X$ be an AM-algebra and let $\alpha >0$. A function
    $f\colon K_X\to \R$ is said to be \emph{$\alpha
    $-homogeneous} if $f(\lambda x^{*})=\lambda ^{\alpha
    }f(x^{*})$ for all $x^{*}\in K_X$ and $\lambda \in \R_+$
    with $\lambda x^{*}\in K_X$. The space
    $C_{ph}^{(\alpha )}(K_X)$ of $\alpha $-homogeneous
    continuous functions on $K_X$ is an AM-space. Moreover, it
    is an AM-algebra with the product given by
    \[
    P(f,g)(x^{*})=\frac{1}{\|x^{*}\|^{\alpha }}
    f(x^{*})g(x^{*})\quad\text{for all }x^{*}\in
    K_X\setminus\{0\}.
    \]
    Indeed, since the norm is weak$^{*}$ continuous on
    $K_X^{*}$, $P$ is a well-defined product on
    $C_{ph}^{(\alpha )}(K_X)$. Moreover, the restriction map
    \[
    \begin{array}{cccc}
    R\colon& C_{ph}^{(\alpha )}(K_X) & \longrightarrow &
    C_b(K_X^{1}) \\
            & f & \longmapsto & f|_{K_X^{1}} \\
    \end{array}
    \]
    defines an \lalg\ isometric embedding. Since
    $C_b(K_X^{1})$ is \lalg\ isometric to the space of
    continuous functions $C(\beta K_X^{1})$, $P$ is the
    AM-algebra product on $X$ (see \cref{rem:pointwise}).
\end{example}

The following are some additional tools to check that a product on an
AM-space is the AM-algebra product.

\begin{lem}
    Let $X$ be an AM-algebra.
    \begin{enumerate}
        \item If $Y$ is a closed subalgebra of $X$, then $Y$ is also
            an AM-algebra (in particular, it is automatically a
            vector sublattice).
        \item Suppose that $T\colon X\to C(K)$ is a lattice isometric
            embedding and that $T(X)$ is a subalgebra. Then $T$ is an
            algebra homomorphism.
    \end{enumerate}
\end{lem}
\begin{proof}~
    \begin{enumerate}
        \item Since $X$ embeds isometrically as a subalgebra in some
            $C(K)$, so does $Y$. In particular, $Y$ is also a
            vector sublattice of $X$, and it is an AM-algebra.
        \item The product induced by $C(K)$ in $T(X)$ is precisely the
            AM-algebra product by \cref{rem:pointwise}.\qedhere
    \end{enumerate}
\end{proof}

\section{AL and AM-spaces with a unique
\texorpdfstring{$f$-algebra}{f-algebra} product}\label{sec:only_zero_prod}

We have first described all the \falg\ products on
AM-spaces, and then studied a particular product (the AM-algebra
product). The existence of this product was characterized by the fact
that the space $W_X$ was maximal. We now turn our heads and look at
the other end. Can $W_X$ be minimal? In other words: can an AM-space
admit a single \falg\ structure? (Namely, the one given by the zero
product.) The following result answers this question in a broader
context. Recall that a Banach lattice $X$ is said to have
\emph{trivial centre} if $\mathcal{Z}(X)$ consists
only of the scalar multiples of the identity.

\begin{prop}
    A Banach lattice with trivial centre and dimension greater than $1$ has a unique \falg\ product.
\end{prop}
\begin{proof}
    Let $X$ be a Banach lattice with trivial centre,
    and let $P$ be an \falg\ product on $X$. We are going to show
    that if $P$ is not the trivial product, then $X$ must be
    one-dimensional. For every $x \in X_+$,
    $P(x,.)$ is an orthomorphism, and therefore a central operator
    \cite[Theorem 3.29]{abramovich_aliprantis2002}. Since the centre is trivial, there
    exists a scalar $\lambda (x)\in \R_+$ such that $P(x,.)=\lambda
    (x)I$. The annihilator of $P$, defined as the set
    \[
    N=\{\, x \in X : P(x,y)=0\text{ for all }y \in X \, \},
    \]
    is a band \cite[Section 3]{huijsmans1991}. If $P$ is not the zero product,
    then $N\neq X$, and therefore its disjoint complement $N^{d}$ is not zero.
    Let $x,y \in (N^{d})_+$. Since $\lambda (x),\lambda (y)\neq 0$ and
    \[
    \lambda (x)y=xy=yx=\lambda (y)x
    \]
    it follows that $N^{d}$ has dimension one. Hence
    $N^{d}$ is a non-trivial projection band. Since the centre of $X$
    is trivial, it must be $X=N^{d}$, and thus $X$ is one-dimensional.
\end{proof}

In particular, every AM-space $X$ of dimension greater than $1$ with trivial centre (see
\cite{goullet_de_rugy1972} and \cite{wickstead1988} for
examples of such spaces) satisfies $W_X=\{0\}$. The converse, however,
is false.

\begin{example}
    Let $X_1$ and $X_2$ be Banach lattices that only admit the trivial
    \falg\ product (for example, $X_1$ and $X_2$
    can be Banach lattices with trivial centre and dimension greater
    than $1$). Their direct
    sum $X=X_1\oplus_\infty X_2$ has non-trivial centre, since both $X_1$ and $X_2$ are non-trivial
    projection bands in $X$. Let $P\colon X\times X\to X$ be an \falg\
    product on $X$. The \falg\ property implies that projection bands
    are closed under multiplication. Hence the restriction
    \[
    P|_{X_1\times X_1}\colon X_1\times X_1\to X_1
    \]
    defines an \falg\ product on $X_1$. By
    assumption, $P|_{X_1\times X_1}=0$.
    Similarly, $P|_{X_2\times X_2}=0$. Let $(x_1,x_2),(y_1,y_2)\in X$
    and compute
    \begin{align*}
        P((x_1,x_2),(y_1,y_2))&=P((x_1,0),(y_1,0))+P((x_1,0),(0,y_2))\\
                              &\phantom{=}+P((0,x_2),(y_1,0))+P((0,x_2),(0,y_2))\\
                              &=P((x_1,0),(0,y_2))+P((0,x_2),(y_1,0))\\
                              &=0,
    \end{align*}
    where in the last equality we are using that the product of
    disjoint elements is $0$. Hence $P=0$.
\end{example}

This prompts a natural question: is there a simple condition
characterizing the AM-spaces with a single \falg\ product? As it turns
out, for AL-spaces there is such a simple characterization. This
characterization will be a simple consequence of the next theorem,
which describes all \falg\ products on AL-spaces.

Before presenting the result, let us recall some notions about atoms
in Banach lattices (see, for instance, \cite[Section 26]{luxemburg_zaanen1971}). An \emph{atom} in a Banach
lattice $X$ is a positive element $a \in X_+$ for which $0\le x\le a$ implies that
$x$ is a scalar multiple of $a$. In this case, the linear span of $a$ is a
projection band. Denote by $P_a$ the associated band projection.
Thus
$P_ax=\lambda _a(x) a$, for a certain $\lambda _a(x)\in \R$. The map
$x\mapsto \lambda _a(x)$ defines a positive linear functional on $X$,
called the \emph{coordinate functional} of $a$. Let $A$ be a maximal
family of norm-one atoms in $X$, and let $B=B(A)$ denote the
band generated by these atoms. When $X$ is Dedekind complete, the
band projection onto $B$ is given by $Px=\bigvee_{a \in A} \lambda
_a(x)a$ for every $x \in X_+$. When $X$ is order continuous,
\[
    Px=\sum_{a \in A}^{} \lambda _a(x)a=\sup_{F \in \Fin(A)}\sum_{a \in
    F}^{}\lambda _a(x)a,
\]
where $\Fin(A)$ denotes the collection of finite subsets of $A$.

Finally, recall that every AL-space is lattice isometric to $L_1(\mu )$,
for some measure $\mu $ \cite[Theorem 7]{kakutani1941AL} (see also
\cite[Theorem 1.b.2]{lindenstrauss_tzafriri1979}).

\begin{thm}\label{thm:ALproducts}
    Let $X$ be an AL-space and let $\{e_i\}_{i \in I}$ be a maximal
    family of norm-one atoms in $X$. Let $e_i^{*}$ be the coordinate
    functional of $e_i$. A product $P\colon
    X\times X\to X$ is an \falg\ product if and only if there exists a
    weight $w \in \el \infty (I)_+$ such that
    \[
    P(f,g)=\sum_{i \in I}^{}w(i)e_i^{*}(f)e_i^{*}(g)e_i\quad\text{for
    all }f,g \in X.
    \]
    This product is submultiplicative if and only if $\|w\|_\infty \le
    1$.
\end{thm}
\begin{proof}
    Identify $X$ with $L_1(\mu )$, for some measure $\mu $.
    Let $P\colon L_1(\mu )\times L_1(\mu )\to L_1(\mu )$ be an \falg\
    product. First assume that $\mu $ is
    $\sigma $-finite. Actually, we may assume that $\mu $ is a
    probability measure, since
    $L_1(\mu )$, with a $\sigma $-finite $\mu $, is lattice isometric to $L_1(\mu ')$ for some
    probability measure $\mu '$.

    For every $f \in L_1(\mu )$, the operator $P(f,.)$ is in the
    centre of $L_1(\mu )$. By \cite[Theorem
    3.34]{abramovich_aliprantis2002}, there exists a unique $w(f) \in L_\infty (\mu
    )$ such that $P(f,g)=w(f) g$ for all $g \in L_1(\mu )$; in fact,
    $w(f) = P(f,\mathbb{1})$. Thus $P(f,g)=P(f,\mathbb{1})g$ and,
    repeating the previous argument with $P(.,\mathbb{1})$, it follows
    that $w(f)=P(f,\mathbb{1})=P(\mathbb{1},\mathbb{1})f$. Therefore
    \begin{equation}\label{eq:L1_prod}
        P(f,g)=P(\mathbb{1},\mathbb{1})fg\quad\text{for all }f,g \in L_1(\mu ).
    \end{equation}

    The map $w\colon L_1(\mu )\to L_\infty (\mu )$, defined by
    $w(f)=P(\one,\one)f$, is a lattice
    homomorphism. Decompose
    $L_1(\mu )=L_1(\nu _1)\oplus L_1(\nu _2)$, where $\nu _1$ is a
    finite atomic measure and $\nu _2$ is a finite
    non-atomic measure. Both $L_1(\nu _1)$ and $L_1(\nu _2)$ are
    projection bands in $L_1(\mu )$. Therefore
    $w|_{L_1(\nu _1)}\colon L_1(\nu _1)\to L_\infty (\mu )$ is also a lattice
    homomorphism. Since the dual $L_\infty (\nu _1)$ of $L_1(\nu _1)$ has no atoms (because
    $\nu _1$ is non-atomic), $L_1(\nu _1)$ does not have non-zero functionals that
    are lattice homomorphisms. But $L_\infty (\mu )$ certainly does.
    So it must be $w|_{L_1(\nu _1)}=0$, that is, $P(\one, \one) \in
    L_1(\nu _2)$. After decomposing
    $f=f_d+f_c, g=g_d+g_c \in L_1(\nu _1)\oplus L_1(\nu _2)$, and
    using that the pointwise product of disjoint functions is zero,
    equation \eqref{eq:L1_prod} ends up being
    \[
    P(f,g)=P(\one, \one)f_dg_d=P(f_d,g_d).
    \]

    Now suppose that $\mu $ is a general measure. For every $f,g \in L_1(\mu  )$, the band generated by $f$ and $g$
    is isometric to $L_1(\mu  ')$, where $\mu  '$ is the restriction of
    $\mu $ to the union of the supports of $f$ and $g$; in particular,
    $\mu '$ is $\sigma $-finite. Since bands are subalgebras, $L_1(\mu
    ')$ is also an \falg\ under the
    inherited product. By the previous argument, $P(f,g)=P(f_d,g_d)$
    in $L_1(\mu ')$, where $f_d$ and $g_d$ denote the projections of
    $f$ and $g$ onto the band generated by the atoms in
    $L_1(\mu ')$. Since $L_1(\mu ')$ is a band in $L_1(\mu )$, $f_d$ and $g_d$ are also the
    projections of $f$ and $g$ onto the band generated by the atoms in
    $L_1(\mu )$.

    As in the statement of the theorem, let $\{e_i\}_{i \in I}$ be a
    maximal family of norm-one atoms in $L_1(\mu )$, and let
    $\{e_i^{*}\}_{i \in I}$ be the associated coordinate functionals.
    If $e_j^{*}(P(e_i,e_i))\neq 0$ were true for some $i\neq j$, then
    \[
    P(e_i,e_i)\wedge e_j\ge e_j^{*}(P(e_i,e_i)) e_j \neq 0
    \]
    contradicting the \falg\ property. So $P(e_i,e_i)$ must be a
    scalar multiple of $e_i$. Let $w \in \el \infty (I)_+$ be defined
    by $w(i)e_i=P(e_i,e_i)$. Note
    that $w$ is indeed bounded, because $P$ is positive on both
    coordinates, and therefore there exists a constant $\|P\|$
    satisfying $\|P(f,g)\|\le \|P\| \|f\| \|g\|$. Using that disjoint
    elements have product zero, we get
    \begin{align*}
        P(f,g)&=P(f_d,g_d)\\&=P\bigg(\sum_{i \in I}^{}e_i^{*}(f)e_i, \sum_{j \in
        I}^{}e_j^{*}(g)e_j\bigg)\\&=\sum_{i \in I}^{}
        e_i^{*}(f)e_i^{*}(g)P(e_i,e_i)\\&=\sum_{i \in
    I}^{}w(i)e_i^{*}(f)e_i^{*}(g)e_i.
    \end{align*}
    Using this expression, it is clear that $P$ is submuliplicative if
    and only if ${\|w\|_\infty \le 1}$.
\end{proof}

\begin{cor}\label{cor:ALspaces_zero}
    Let $X$ be an AL-space. The only possible \falg\ product on $X$ is
    the zero product if and only if $X$ has no atoms.
\end{cor}

For a general Banach lattice $X$, it is true that if $X$ has an atom,
then there exists a non-trivial \falg\ product on $X$.
This is proved after the following lemma, which may be of
interest by itself.

\begin{lem}\label{lem:non-trivial_local}
    Let $X$ be a vector lattice. If a projection band $B$ in $X$
    admits a non-trivial \falg\ product, so does $X$.
\end{lem}
\begin{proof}
    Suppose $B$ is a projection band in $X$ with a non-trivial
    \falg\ product $P$. Let $Q\colon X\to B$ be the band projection
    associated with $B$, and define $\tilde P\colon X\times X\to X$ by
    $\tilde P(x,y)=P(Qx,Qy)$. This is certainly a product on $X$. We
    will be done once we show that $\tilde P$ satisfies the \falg\ property. If $x\wedge y=0$
    in $X$, then $Qx \wedge Qy=0$, because $Q$ is a lattice
    homomorphism. Let $z \in X_+$. Since $P$ is an \falg\ product,
    $P(Qz,Qx)\wedge Qy=0$. But also $P(Qz,Qx)\wedge (y-Qy)=0$, because
    $P(Qz,Qx) \in B$ and $y-Qy \in B^{d}$. In conclusion,
    \[
    0\le \tilde P(z,x)\wedge y\le P(Qz,Qx)\wedge Qy+P(Qz,Qx)\wedge (y-Qy)=0,
    \]
    so $\tilde P$ does satisfy the \falg\ property.
\end{proof}

\begin{prop}
    If a vector lattice has an atom, then it has a non-trivial \falg\
    product.
\end{prop}
\begin{proof}
    Let $X$ be a vector lattice and let $a \in X_+$ be an atom. Then
    $\barespn \{a\}$ is a projection band that admits a non-trivial \falg\
    product (for instance, $P(\lambda a,\mu a)=(\lambda \mu )a$ for all
    $\lambda ,\mu \in \R$). By \cref{lem:non-trivial_local}, $X$
    admits a non-trivial \falg\ product.
\end{proof}

As we already know, the converse of previous result is true for
AL-spaces, but not in general (it fails already in $C[0,1]$). In the
end, we were not able to find an analogue of \cref{cor:ALspaces_zero}
for AM-spaces.

\begin{question}
    Is there a simple characterization of the AM-spaces $X$
    for which the only \falg\ product is the zero product (i.e., for
    which $W_X=\{0\}$)?
\end{question}

\section*{Acknowledgements}

I would like to thank Pedro Tradacete for valuable
discussions and corrections. I would also like to thank Marcel de Jeu
for his helpful comments and for pointing out item \textit{(vii)}
in \cref{thm:charAMalg}. Finally, I would like to thank the anonymous
referee for suggestions that significantly improved the clarity of the
paper, and for pointing out the elegant proof of
\cref{prop:unique_general}.

\emergencystretch=1em
\printbibliography

\end{document}